\documentclass[11pt,a4paper,twoside]{amsart}
\usepackage{amsfonts, amssymb,indentfirst}

\newtheorem{prop}{Proposition}[section]
\newtheorem{theorem}[prop]{Theorem}
\newtheorem{lemma}[prop]{Lemma}

\newtheorem{coro}[prop]{Corollary}
\newtheorem{remark}[prop]{Remark}
\newenvironment{rem}{\begin{remark}\rm}{\end{remark}}
\newcommand{\cqd}{\hfill$\Box$}
\newcommand{\mor}[0]{\operatorname{Mor}}
\newcommand{\quot}[0]{\operatorname{Quot}}
\renewcommand{\deg}[0]{{\it deg}}
\renewcommand{\dim}[0]{\operatorname{dim}}
\newcommand{\codim}[0]{\operatorname{codim}}
\newcommand{\im}[0]{\operatorname{Im}}

\title[RATIONAL RULED SURFACES AND GROMOV-WITTEN INVARIANTS] {THE
  DEGREE OF THE VARIETY OF RATIONAL RULED SURFACES AND GROMOV-WITTEN
  INVARIANTS} \author[Cristina Mart{\'\i}nez]{Cristina Mart{\'\i}nez}

\address{Departamento de Matematicas, Facultad
de Ciencias, Universidad Autonoma de Madrid, Madrid, 28049, Spain}
 \email{cristina.martinez@uam.es}
\begin{document}
\maketitle

\begin{abstract}
We compute the degree of the variety parametrizing rational ruled
surfaces of degree $d$ in $\mathbb{P}^{3}$ by relating the problem
to Gromov-Witten invariants and Quantum cohomology.
\end{abstract}


\begin{center}
\section{Introduction}
\end{center}
\normalsize \setcounter{equation}{0}

In 1986, D.F. Coray and I. Vainsencher computed the degree of
certain strata of the variety parametrizing ruled cubic surfaces
\cite{CV}. Here we generalize their result and compute the degree
of the variety parametrizing rational ruled surfaces of degree $d$
in $\mathbb{P}^{3}$ in the projective space parametrizing all
surfaces of degree $d$. We approach this enumerative problem by
fixing a suitable parameter space for the objects that we want to
count, and expressing the locus of objects satisfying given
conditions as a zero-cycle on this parameter space. We then need
to evaluate the degree of this zero-dimensional cycle class. This
is possible in principle whenever the Chow group of cycles up to
numerical equivalence of the parameter space is known in terms of
generators and relations.

We use $\mor(\mathbb{P}^{1}, G(2,4))$, the variety of morphisms
from $\mathbb{P}^{1}$ to $G(2,4)$, as a parameter space for
rational ruled surfaces of degree $d$. Here $G(2,4)$ denotes the
Grassmannian of lines in $\mathbb{P}^3$. By the universal property
of the Grassmannian, we can identify a rational ruled surface in
$\mathbb{P}^{3}$ with a rational curve in $G(2,4)$. More
precisely, the points of this scheme are parametrized ruled
surfaces. The objects we intend to count are the images of these
morphisms, which are rational curves of fixed degree $d$ in the
Grassmannian $G(2,4)$. To eliminate the data of the
parametrization due to the action, of the group PGL($1$) on
$\mathbb{P}^{1}$, we impose three point conditions.

The variety of morphisms of a fixed degree from $\mathbb{P}^{1}$
to the Grassmannian $G(2,4)$ is not a compact parameter space. The
Grothendieck Quot scheme parametrizing rank 2 and degree $d$
quotients of a trivial vector bundle
$\mathcal{O}^{4}_{\mathbb{P}^{1}}$ is a compactification of this
variety as shown in \cite{Str}. This Quot scheme has very good
properties, it is a fine moduli space equipped with a universal
element, it is a smooth, irreducible scheme. Unfortunately, some
of the divisors we want to intersect have a common component
contained in the boundary of the Quot scheme. Therefore, our
intersection problem does not have enumerative meaning on the Quot
scheme.

We instead use the coarse moduli space
$\overline{M}_{0,n}(G(2,4),d)$ of Kontsevich stable maps from
$n-$pointed genus 0 curves to the Grassmannian $G(2,4)$
representing $d$ times the positive generator of the homology
group $H_{2}(G(2,4),\mathbb{Z})$. In the Kontsevich space our
classes intersect properly and have enumerative meaning.  Since
the Grassmannian $G(2,4)$ is a homogenous variety transversality
arguments imply a relationship between Gromov-Witten invariants
and enumerative geometry. We proceed in some sense in an opposite
way to Bertram in \cite{Ber1} where he defines the Gromov-Witten
invariants and then he reinterprets them as intersections of
generalized Schubert cohomology classes in a Grothendieck Quot
scheme.

The quantum cohomology ring of a variety is defined in terms of
intersection data (the Gromov-Witten invariants) on the spaces of
holomorphic maps from pointed curves of genus zero to the variety.
The quantum cohomology of $G(2,4)$, $QH^{*}(G(2,4))$, is an
algebra over a polynomial ring in one variable. It has been
described in \cite{FI}, where the authors present some of the line
and conic quantum numbers. The solutions to our enumerative
problem occur as coefficients in the multiplication table of
$QH^{*}(G(2,4))$.

The associativity of the quantum cohomology ring gives many relations
among the Gromov-Witten invariants and allows us to determine the
invariants from a few basic ones. In this work, we have applied
Farsta, a computer program due to Andrew Kresch, which computes
quantum numbers using associativity relations, to find the invariants
we are interested in.

\vspace{0.5cm}
\subsubsection*{Acknowledgements}
This work was carried out under the supervision of Profesor Rafael
Her\-n{\'a}n\-dez at Aut{\'o}noma University of Madrid. I would
like to thank Aaron Bertram for some encouraging and clarifying
conversations in Lille, and Paolo Aluffi and Andrew Kresch for
their help with farsta program. I also thank the referee for a
long list of corrections and suggestions and Izzet Coskun and
Davesh Maulik for a careful reading of a previous version of this
paper. My work has been supported by an FPI grant BFM 2000-0024.

\vspace{0.5cm} \small \subsubsection*{\small Notation} $A_{d}X$
and $A^{d}X$ can be taken to be the Chow homology and cohomology
groups for homogeneous varieties. For $\beta\in A_{k}X$,
$\int_{\beta}c$ is the degree of the zero cycle obtained by
evaluating $c_{k}$ on $\beta$, $c\in A^{*}X$. We use cup product
$\cup$ for the product in $A^{*}X$. A closed subvariety $T$ of
$G(2,4)$ of pure codimension $c$ determines classes in
$A_{n-c}(G(2,4))$ and $A^{c}(G(2,4))$ via the duality isomorphism.
Both of these classes are denoted by $[T]$. We denote the
automorphism group of $\mathbb{P}^{n}$ by PGL($n$).

Throughout this paper, we work over the field of complex numbers
$\mathbb{C}$.

\normalsize
\section{{\large}Rational ruled surfaces}
\subsection{Some basic preliminaries}\label{themorph}
\setcounter{equation}{0} A {\it rational ruled surface} is the
projectivization of a locally free sheaf of rank 2 on
$\mathbb{P}^{1}$, $Y=\mathbb{P}(E)$, together with the projection
morphism $\pi:\mathbb{P }(E)\rightarrow \mathbb{P}^{1}$ (see $V.2$
of \cite{Har}). Classically the term referred to a birational
image of $Y$ in $\mathbb{P}^{3}$ where the fibers of $\pi$ are
mapped to lines. We will describe how such an $f$ can be obtained.

If $E$ is isomorphic to the locally free sheaf
$\mathcal{O}_{\mathbb{P}^{1}}^{2}$ in the open set $U\subset
\mathbb{P}^{1}$, then
$$\pi^{-1}U= \mathbb{P}(E)|_{U}\cong \mathbb{P}(\mathcal{O}^{2}_{\mathbb{P}^{1}})
\mid_{U}\cong\mathbb{P}^{1}\times U.$$ 

This means that $\pi:\mathbb{P}(E)\rightarrow
 \mathbb{P}^{1}$ is a locally trivial fibration. By a Theorem
 of Grothendieck every locally free sheaf over $\mathbb{P}^{1}$
decomposes as a direct sum of linear sheaves.

{\it The degree of $E$} is the degree of the invertible sheaf
$\bigwedge^{2}E$. If $E\cong \mathcal{O}_{\mathbb{P}^{1}}(a)\oplus
\mathcal{O}_{\mathbb{P}^{1}}(b)$ then
$$\bigwedge^{2} E \cong \mathcal{O}_{\mathbb{P}^{1}}(a)
\otimes \mathcal{O}_{\mathbb{P}^{1}}(b)\cong
\mathcal{O}_{\mathbb{P}^{1}}(a+b),\ \ d:=a+b. $$

We can suppose $a,b\geq 0$ (see V.2.2. of \cite{Har}). In
particular this implies $E$ is generated by global sections, or
equivalently $E$ is given by a quotient of
$\mathcal{O}^{n+1}_{\mathbb{P}^{1}}$:
\begin{equation}\label{eq}
0\rightarrow \mathcal{N} \rightarrow
\mathcal{O}^{n+1}_{\mathbb{P}^{1}}\rightarrow E\rightarrow 0
\end{equation}
Moreover, we will suppose $a,b >0$, therefore we are excluding the
cone of degree $d$, \  $\mathbb{P}(\mathcal{O}_{\mathbb{P}^{1}}\oplus
\mathcal{O}_{\mathbb{P}^{1}}(d))$.

\noindent The sequence (\ref{eq}) induces the morphism
$\mathbb{P}(E)\stackrel{i}{\rightarrow}\mathbb{P}(\mathcal{O}^{n+1}_{\mathbb{P}^{1}})$.
But $\mathbb{P}(\mathcal{O}^{n+1}_{\mathbb{P}^{1}})\cong
\mathbb{P}^{n}\times \mathbb{P}^{1}$ and by composing this
morphism with the projection to the first factor, we
obtain a morphism
\begin{equation} \label{eq3}
f_{n}:\mathbb{P}(E)\rightarrow \mathbb{P}^{n}.
\end{equation}
This is a projective morphism by definition, but in general it is
not a birational map. We denote by $X$ the image of $f_{n}$.

\noindent It is easy to see the morphism $f_{n}$ is associated to
the linear sheaf $\mathcal{O}_{\mathbb{P}(E)}(1)$ and maps fibers
of $\pi$ into lines. Moreover, the two projections
\begin{equation}\pi_{a}:\mathcal{O}_{\mathbb{P}^{1}}(a)\oplus
\mathcal{O}_{\mathbb{P}^{1}}(b)\rightarrow
\mathcal{O}_{\mathbb{P}^{1}}(a)\end{equation}
\begin{equation}\pi_{b}:\mathcal{O}_{\mathbb{P}^{1}}(a)\oplus
\mathcal{O}_{\mathbb{P}^{1}}(b)\rightarrow
\mathcal{O}_{\mathbb{P}^{1}}(b)\end{equation} correspond to
sections (II.7.12. of \cite{Har}) which are mapped to curves
$C^{a}$ and $C^{b}$ of degrees $a$ and $b$ respectively. $X$ can
be obtained from $C^{a}$ and $C^{b}$: for each point $p\in C^{a}$
there is a unique point $q\in C^{b}$ such that $\pi(p)=\pi(q)$.
Then $X$ is the union over all $p\in C^{a}$ of the lines
$\overline{pq}$, (see \S 5.6 of \cite{Fri}).

\noindent For $a=b=1$, there are two distinct families of lines
and $f_{n}(\mathbb{P}(\mathcal{O}_{\mathbb{P}^{1}}(1)\oplus
\mathcal{O}_{\mathbb{P}^{1}}(1)))$ is the quadric $Q$ in
$\mathbb{P}^{n}$. In this case, there are two different morphisms
$Q\rightarrow \mathbb{P}^{1}$ corresponding to the two
projections.

\begin{lemma}
The morphism (\ref{eq3}) is a finite morphism.
\end{lemma}

{\it Proof.} Let us consider the intersection number
\[f^{*}\mathcal{O}_{\mathbb{P}^{n}}(1)\cdot f^{*}\mathcal{O}_{\mathbb{P}^{n}}(1)=
\mathcal{O}_{\mathbb{P}(E)}(1)^{2}.\] Let  $\xi$ denote
 the divisor class corresponding to $\mathcal{O}_{\mathbb{P}(E)}(1)$
 in the Chow ring $A(\mathbb{P}(E))$.  By the definition of the
first Chern class $c_{1}(E)\in
A^{1}(\mathbb{P}(E))$, we have \[\pi^{*}c_{0}(E)\xi^{2}-
\pi^{*}c_{1}(E)\xi=0,\] and $c_{0}(E)=1$. Hence,
\begin{displaymath}
\xi^{2}=\pi^{*}c_{1}(E)\,\xi =deg\, E=d>0.
\end{displaymath}

Now from the fact that $\dim X=\dim \mathbb{P}(E)=2$, it follows
that $f_{n}$ is generically finite by the fiber-dimension theorem
(I. 6.3.7, \cite{Sha}). For each $y\in \mathbb{P}^{n}$, the
pre-image $\pi_{2}^{-1}(y)$ is a section of
$\mathbb{P}(\mathcal{O}^{n+1}_{\mathbb{P}^{1}})\rightarrow
\mathbb{P}^{1}$, therefore a quotient
$\mathcal{O}^{n+1}_{\mathbb{P}^{1}}\twoheadrightarrow
\mathcal{O}_{\mathbb{P}^{1}} $. $\pi_{2}^{-1}(y)$ is either a
0-dimensional fiber or a section of $\mathbb{P}(E)$. If it were a
section, we would have,
$$ \begin{picture}(100,60)
\put(0,50){$\mathcal{O}_{\mathbb{P}^{1}}^{n+1}\twoheadrightarrow
\mathcal{O}_{\mathbb{P}^{1}}$} \put(10,40){\vector(1,-1){29}}
\put(50,10){\vector(0,2){35}}\put(45,0){$E$}
\end{picture} $$

Since $Hom\,(E,\mathcal{O})=H^{0}(E^{\vee})=0$ and $a,b>0$,
$\pi_{2}^{-1}(y)\subset \mathbb{P}(E)$ cannot be a section.
Therefore, $\pi_{2}^{-1}(y)$ is a finite set of points.  This
shows that the morphism is quasi-finite. Actually, it is finite
(see exercise 3.11.2 of \cite{Har}). We also conclude the
following formula connecting the degree of the morphism with the
degree of $X$,
\begin{equation}\label{degree}
\deg(X)\cdot \deg(f)=f^{*}(\mathcal{O}_{\mathbb{P}^{n}}(1))^{2},
\end{equation}
where $\deg(f)$, is the degree of the field extension
$K(f(X))\subset K(\mathbb{P}(E))$.\cqd

\begin{rem} When $n=h^{\circ}(\mathbb{P}^{1},E)=d+2$, the surface is said to be \lq\lq
  linearly normal". The surfaces obtained as the images of the above
  morphisms are all projections of the linearly normal surface. In
  this particular case, $deg\,X=\mathcal{O}_{\mathbb{P}(E)}(1)^{2}=d$
  (see \S 1 of \cite{EH}). Therefore formula (\ref{degree}) implies
  that $\deg(f)=1$. Hence, $f_{d+1}:\mathbb{P}(E)\rightarrow
  \mathbb{P}^{d+1}$ is generically injective, that is, birational onto
  its image.
\end{rem}

\subsection{A parameter space for rational ruled surfaces of degree
$d$ in $\mathbb{P}^{3}$.}In order to solve a specific counting
problem, one of the most successful approaches is to apply the
methods of intersection theory to parameter spaces. We will put
the family of rational ruled surfaces in one to one correspondence
with the points of an algebraic variety, the parameter space or
moduli space.

\begin{lemma}\label{bijection}
There exists a bijective correspondence \mbox{between} the sets of
isomorphism classes of rank 2 and degree $d$ quotients on
$\mathbb{P}^{1}$ of $\mathcal{O}^{n+1}_{\mathbb{P}^{1}}$ and the
set of morphisms of degree $d$ of $\mathbb{P}^{1}$ in the
Grassmannian of 2-dimensional subspaces of a
$\mathbb{C}-$vector space of dimension $n+2$,
$\mor(\mathbb{P}^{1},Gr(2,n+1))$.
\end{lemma}
{\it Proof.}
 Let $E$ be a bundle of rank 2 on $\mathbb{P}^{1}$, with
$\mathcal{O}^{n+1}_{\mathbb{P}^{1}}\rightarrow E\rightarrow 0$,
and let us consider the universal exact sequence on the
Grassmannian,
\begin{equation} \label{eq6}
0\rightarrow \mathcal{S}\rightarrow
\mathcal{O}^{n+1}_{G}\rightarrow \mathcal{Q}\rightarrow 0.
\end{equation}
Then by the universal property of the Grassmannian, there exists a
unique morphism $\varphi:\mathbb{P}^{1}\rightarrow G(2,n)$, so
that
$$
(\varphi^{*}\mathcal{O}^{n+1}_{G}\rightarrow
\varphi^{*}\mathcal{Q}\rightarrow
0)=\mathcal{O}^{n+1}_{\mathbb{P}^{1}}\rightarrow
E_{\mathbb{P}^{1}}\rightarrow 0.
$$

\noindent The degree of $\varphi$ is defined as, $$ \deg
(\varphi):=\deg(\varphi^{*}\mathcal{Q}).$$

\noindent The pull-back under $\varphi$ of the universal quotient
on the Grassmannian $\varphi^{*}\mathcal{Q}$, is a bundle of rank 2 on
$\mathbb{P}^{1}$ locally isomorphic to
$\mathcal{O}_{\mathbb{P}^{1}}(a)\oplus
\mathcal{O}_{\mathbb{P}^{1}}(d-a)$ where $(0\leq a\leq
[\frac{d}{2}])$. On the other hand
$c_{1}(\mathcal{Q})=\mathcal{O}_{G}(1)$ and
$\varphi^{*}(c_{1}(\mathcal{Q}))=c_{1}(\varphi^{*}(\mathcal{Q}))=c_{1}(E),$
its first Chern class, therefore $\deg(\varphi)=deg(E)=d$. \cqd

\vspace{0.5cm} Now we concentrate our attention on the case $n=3$,
i.e. quotients
\begin{equation}\label{cociente}
\mathcal{O}^{4}_{\mathbb{P}^{1}}\rightarrow E\rightarrow 0.
\end{equation} In this case the morphism $f_{3}$ given in
\ref{themorph}, maps the rational ruled surface $\mathbb{P}(E)$
into $\mathbb{P}^{3}$ and, by the previous identification, we can
see it as a curve of degree $d$ in $G(2,4)$.

If we fix the degree, $d$, and the rank, 2, of a locally free
sheaf $E$ on $\mathbb{P}^{1}$, we are fixing its Hilbert
polynomial, $$\label{eq10} P(t)=\chi(E(t))=2t+d+2.
$$
The moduli, $\quot(\mathbb{P}^{1},P(t))$, of quotients with fixed
Hilbert polynomial $P(t)$ is a fine moduli space by a theorem of
Grothendieck \cite{Gro}. We will denote it as $R_{d}$. We observe
that the quotient $\mathcal{O}^{4}_{\mathbb{P}^{1}}\rightarrow
E\rightarrow 0$, determines a point $q\in
\quot(\mathbb{P}^{1},P(t))$ and a morphism
$f_{q}:\mathbb{P}^{1}\rightarrow G(2,4)$ by the universal property
of the Grassmannian. By definition, there is a universal quotient,
\begin{equation}\label{eq11}
\mathcal{O}^{4}_{\quot\times\mathbb{P}^{1}}\rightarrow
\mathcal{E}_{\quot\times \mathbb{P}^{1}},
\end{equation}
with the property that for all $k-$schemes $S$, the set of
morphisms $f:S\rightarrow \quot$ is in one to one correspondence
with the set of isomorphism classes of short exact sequences over
$S\times \mathbb{P}^{1}$,
$$\mathcal{O}^{4}_{S\times \mathbb{P}^{1}}\rightarrow
\mathcal{E}_{S\times \mathbb{P}^{1}}\rightarrow 0 $$  where $\
\mathcal{E}_{S\times \mathbb{P}^{1}}$  is flat over $S$  with
Hilbert\ polynomial
$$\chi(\mathcal{E}_{\{ s\}})=2t+d+2 \rm \ on \ the \ fibers \ of\
\pi_{s}:S \times \mathbb{P}^{1}\rightarrow S.$$

Equivalently, $\mathcal{E}_{\{s\}\times \mathbb{P}^{1}}$ has rank
2 and degree $d$ for every $s\in S$. There exists an open neighborhood
$U$ of
$\{q\}$ on the Quot scheme so that $\mathcal{E}_{U\times
\mathbb{P}^{1}}$ is locally free of rank 2 and degree $d$.

Let $R^{0}_{d}$ be the maximal open set where $\mathcal{E}$ is
locally free and of constant rank 2, then $R^{0}_{d}\cong
\mor_{d}(\mathbb{P}^{1},G(2,4))$. $R_{d}$ is a natural
compactification of $R^{0}_{d}$ as shown in \cite{Str}. In
\cite{Str} it is proved that $R^{0}_{d}$ is a quasi-projective,
smooth, rational and irreducible variety of dimension $4d+4$. This
means that quotients (\ref{cociente}) are points of the scheme
$R^{0}_{d}$. That is, it parametrizes ruled surfaces $X\subset
\mathbb{P}^{3}$ together with a projection map $\pi:X\rightarrow
\mathbb{P}^{1}$ whose fibers are isomorphic to $\mathbb{P}^{1}$.

\begin{lemma} \label{finemod} The space $R^{0}_{d}$ is a fine moduli space for
degree $d$ maps from $\mathbb{P}^{1}$ to $G(2,4)$.
\end{lemma}
{\it Proof.} Let us consider the evaluation map from
$R^{0}_{d}\times \mathbb{P}^{1}$ to the Grassmannian $G(2,4)$, and the
projection map to the first component.

\begin{equation}\label{diagram}\begin{array}{ccccccccc}
 &   & R^{0}_{d}\times \mathbb{P}^{1} & \stackrel{e}{\rightarrow} & G(2,4) & & & & \\
 &   & \pi_{1} \downarrow  & & & & & & \\
  & & R^{0}_{d} & & & & & & \\
\end{array}\end{equation}

\noindent This family is a universal family. By the universal
property of the Grassmannian, the pull back under $e$ of the
universal exact sequence on $G(2,4)$ (\ref{eq6}), gives us a universal
exact sequence on $R^{0}_{d}$ which is the restriction of the
universal one on $R_{d}\times \mathbb{P}^{1}$.  \cqd

\vspace{0.5cm} We will denote by $R^{00}_{d}$ the Zariski open set
in $R^{0}_{d}$ corresponding to a morphism birational onto its
image. This open set is precisely the set of automorphism free
maps.

\begin{lemma}\label{codim} If $d\geq 2$,  the complement
of $R^{00}_{d}$ in $R^0_d$, the locus parametrizing multiple covers, has
codimension at least 2.
\end{lemma}

{\it Proof.} A birational morphism on $R^{0}_{d}$ corresponds to a
birational map to $\mathbb{P}^{3}$, (see \ref{themorph} and
\ref{bijection}), and consequently to a rational ruled surface of
degree $d$ in $\mathbb{P}^{3}$. For $d\geq 2$, formula
(\ref{degree}) implies that the maps in the complement of
$R^{00}_{d}$ are the multiple covers.
 For $k|d$, every $k-$multiple cover factorizes as,
$$
\mathbb{P}^{1}\stackrel{\rho}{\rightarrow}\mathbb{P}^{1}\stackrel{\xi}
{\rightarrow} G(2,4),
$$
where $\rho$ is a $k-$sheeted cover of $\mathbb{P}^{1}$ and
$\xi\in \mor_{d/k}(\mathbb{P}^{1},G(2,4))$. Then there is a
natural morphism,
$$
\begin{array}{cccc}
\mor_{k}(\mathbb{P}^{1},\mathbb{P}^1) \ \times &
\mor_{d/k}(\mathbb{P}^{1},G(2,4)) & \rightarrow &
\mor_{d}(\mathbb{P}^{1},G(2,4)) \\ \rho & \xi & & \rho\circ \xi
\end{array}
$$
The dimension of the product is $(2k+1)+(4\frac{d}{k}+4)+3$. The
fibers of this morphism, PGL(1), are of dimension 3, therefore the
image has codimension,
$$4d+4+3-(2k+1)-(4\frac{d}{k}+4)=4d-4\frac{d}{k}+2-2k.$$
The study of this function for $2\leq k\leq d$ shows that the
maximum value must be attained at the end points $k=2,d$. If
$k=2$,
$$4d-2d-4=2\,(d-1)\geq 0.$$ If $k=d$ since $d\geq 2$, $$2d-4=2\,(d-2)\geq 0.$$
This shows that the complement $(R^{00}_{d})^{c}$ is of
codimension at least 2. \cqd

\begin{rem}
In addition to the problem of multiple covers, reparametrizations
of the same curve in $G(2,4)$ are considered distinct objects in
$R^{0}_{d}$. This means, that the actual space for rational ruled
surfaces ought to be $R^{00}_{d}$/PGL(1), and this is in fact a
quotient in the sense of Mumford.
\end{rem}

\begin{prop}\label{morphism}
There exists a morphism from the variety $R^{00}_{d}$ to the
Hilbert scheme $\mathbb{P}^{{d+3 \choose 3}-1}$ of surfaces in
$\mathbb{P}^{3}$ of degree $d$.
\end{prop}
{\it Proof.} From the universal quotient (\ref{eq11}) we obtain a
morphism
$$
\mathbb{P}(\mathcal{E}_{ \mathbb{P}^{1}\times
R^{0}_{d}})\rightarrow
\mathbb{P}(\mathcal{O}^{4}_{\mathbb{P}^{1}\times R^{0}_{d}})\cong
\mathbb{P}^{1}\times R^{0}_{d}\times \mathbb{P}^{3} .$$ Projecting
to the  last two components, we obtain a morphism
$$\label{fami-uni} \mathbb{P}(\mathcal{E}_{R^{0}_{d}\times
\mathbb{P}^{1}})\rightarrow  R^{0}_{d}\times \mathbb{P}^{3}.
$$

\noindent For all $q\in R^{0}_{d}$, we have that
$\mathbb{P}(E_{q})\rightarrow \{q\}\times \mathbb{P}^{3}$ with
$E_{q}:=\mathcal{E}_{q\times \mathbb{P}^{1}}$, is the rational
surface corresponding to that point by the morphism constructed in
\ref{themorph}. Although we have seen this morphism is not always
birational, the map \ref{fami-uni} restricted to the open set
$R^{00}_{d}\subset R^{0}_{d}$ is a birational map. Now composing
with the projection morphism to the first component, we have a
morphism,

\begin{equation}
\mathbb{P}(\mathcal{E}_{R^{00}_{d}\times
\mathbb{P}^{1}})\stackrel{\varphi}{\rightarrow} R^{00}_{d}.
\end{equation}

\noindent We now consider the image of this morphism,
\begin{equation} \label{flat-fami}
\varphi(\mathbb{P}(\mathcal{E}_{R^{00}_{d}\times
\mathbb{P}^{1}}))\stackrel{i}{\hookrightarrow}
R^{00}_{d}.\end{equation} We observe that for each $r\in
R^{00}_{d}$, $i^{-1}(r)\subset \mathbb{P}^{3}$ is a rational ruled
surface with cons\-tant Hilbert polynomial (see ex.V.1.2 of
\cite{Har}), therefore (\ref{flat-fami}) is a universal family. By
the universal property of the Hilbert scheme parametrizing all
surfaces of degree $d$ in $\mathbb{P}^{3}$, $\mathbb{P}^{{d+3
\choose 3}-1}$, there exists a unique morphism,
\begin{equation}
R^{00}_{d}\stackrel{\phi}{\rightarrow}\mathbb{P}^{{d+3 \choose
3}-1},
\end{equation}
which sends a point $q\in R^{00}_{d}$ to the corresponding surface
$f_{3}(\mathbb{P}(E_{q}))$, where $f_{3}$ is the morphism defined
in \ref{themorph}. \cqd

\section{The degree of the variety of rational ruled surfaces}
\setcounter{equation}{0} \label{degree}

We want to compute the {\bf degree} of the {\bf variety}
parametrizing the family of {\bf rational ruled surfaces} of
degree $d$ in $\mathbb{P}^{3}$
 in the projective space of surfaces in $\mathbb{P}^{3}$ of
 degree $d$, $\mathbb{P}^{{d+3 \choose
 3}-1}$.

 We will express the locus of objects satisfying
 given geometric conditions as a zero-cycle on the
 parameter space we have fixed for rational ruled surfaces. In other
 words, we want to compute {\it the degree $N_{d}$} of
 $\im(\phi)$ for $d\geq 3$ where $\phi$ is the morphism,

\begin{equation}\label{morfismo}
R^{00}_{d}\stackrel{\phi}{\rightarrow}\mathbb{P}^{{d+3 \choose
3}-1},
\end{equation}
defined in \ref{morphism}.

In order to compute this degree, the divisors we have to intersect
are des\-cri\-bed geometrically as sets of rational curves
verifying certain incidence conditions with some Schubert
va\-rie\-ties on the Grassmannian. For this reason, we first
describe the cohomology ring of the Grassmannian. We consider
again the universal sequence on the Grassmannian $G(2,4)$:

\begin{equation}
0\rightarrow \mathcal{S}\rightarrow \mathcal{O}_{G}^{4}\rightarrow
\mathcal{Q}\rightarrow 0.
\end{equation}
We will represent the special Schubert cycles on $G(2,4)$ as
$$
T_{1}=c_{1}(\mathcal{Q})\ \mbox{lines  meeting  a  given  line,}
$$
$$T_{b}=c_{2}(\mathcal{S})\  \mbox{ lines  contained  in  a  given
plane,}$$
$$ T_{a}=c_{2}(\mathcal{Q})\ \mbox{lines containing a
given point.} $$

Also, $T_{3}\in H_{6}(G(2,4),\mathbb{Z})$ will stand for the class
of a
line and $T_{4}$ for the class of a point. 

\noindent We now give a description of the two Weil divisors on
$R^{0}_{d}$ we are interested in. These are constructed by means
of the evaluation map in (\ref{diagram}):
\begin{enumerate}
\item[(A)] The locus of morphisms whose image meets an $a-$plane in
  the Grassmannian
associated to a point $P_{i}\in \mathbb{P}^{3}$. We denote it by $D_{i}$.
$$ \label{div1}
\{\varphi \in R^{0}_{d}|\ \  e(t,\varphi)\cap T_{a_{i}}\neq
\emptyset\}.
$$

\item[(B)] The set of morphisms $Y_{i}$ sending a fixed point
$t_{i}\in \mathbb{P}^{1}$ to a hyperplane $T_{1}$ on the
Grassmannian.
$$\label{div2} \{\varphi \in R^{0}_{d}|\
e(t_{i},\varphi)\in T_{1}\ \mbox{for  a  fixed } \  t_{i}\in \mathbb{P}^{1}
\}.
$$
\end{enumerate}

\begin{theorem}\label{teorema1} For generic choices of
$a-$planes and hyperplanes on the Grassmannian $G(2,4)$, the
degree $N_{d}$ coincides with the intersection of the divisors:
\begin{equation}\label{Nd}
W:= Y_{1}\cap Y_{2}\cap Y_{3}\cap D_{1}\cap \ldots \cap D_{4d+1}
\end{equation}
or equivalently, the product in the Chow ring $A(R_{d}^{0})$ of
the one cycles:
\begin{equation}
[Y_{1}][Y_{2}][Y_{3}][D_{1}]\ldots [D_{4d+1}].
\end{equation}
\end{theorem}

\noindent Although the divisors we intersect are in $R^{0}_{d}$,
we will see in the following lemma the intersection is in fact in
the open set $R^{00}_{d}$ and it is a finite number of reduced
points which coincides with the degree $N_{d}$ of $\im(\phi)$.

\noindent  The Grassmannian of lines $G(2,4)$ is a homogeneous
variety under the action of the group PGL(3).
\begin{lemma}\label{intnum}
The intersection $W$ consists of a finite number of reduced points
supported in the locus $R^{00}_{d}$ of maps without automorphisms.
\end{lemma}
{\it Proof.} Let $G(2,4)^{4d+4}=G(2,4)\times
\stackrel{4d+4}{\ldots}\times G(2,4)$ be the product of $4d+4$
factors equal to $G(2,4)$. We fix three points to be $0, 1,
\infty$ in $\mathbb{P}^{1}$, and we consider the multiple
evaluation map,
\begin{equation}\label{multev} R^{0}_{d}\times\{0\}\times \{1\}\times \{\infty\}
\times \mathbb{P}^{1}\times \stackrel{4d+1}{\ldots}\times
\mathbb{P}^{1}\times \stackrel{\overline{e}}{\rightarrow}
G(2,4)^{4d+4}.\end{equation} Let $\Upsilon$ be the cycle on
$G(2,4)^{4d+4}$ given by the product of the irreducible varieties
on $G(2,4)$ associated to the divisors in (\ref{Nd}), i.e.
$W=\overline{e}\,^{-1}(\Upsilon)$.

\noindent We consider the action of the product of $4d+4$ copies
of PGL(3) on $G(2,4)^{4d+4}$. Let $D_{i,j}$ be the diagonal on
$G(2,4)^{4d+4}$ determined by the factors $i$ and $j$. Restricting
to the complement of $D_{i,j}$, we get a transitive action. We
apply Kleiman's Theorem \cite{Kle} to the complement
$(R^{00}_{d})^{c}$. We have seen in Lemma \ref{codim} that this is
a closed subvariety of codimension at least 2 in $R^{0}_{d}$.
Kleiman's Theorem applied to
\begin{equation}
\begin{array}{ccc}
 & & (R^{00}_{d})^{c}\\
  & & \downarrow \\
  \Upsilon & \hookrightarrow & G(2,4)^{4d+4}\backslash \{D_{i,j}\}
\end{array}
\end{equation}
tells us that the intersection is supported in
$R^{00}_{d}$.

\noindent The intersection $W$ is transverse, smooth and the
corresponding class in the Chow ring $A(R^{0}_{d})$ is a top
degree class. This is a consequence of the usual Bertini Lemma
(see II.8.18 of \cite{Har}).\cqd

\vspace{0.5cm}
 \noindent {\it Proof.} (of Theorem \ref{teorema1})
For a general hyperplane $H$ in $\mathbb{P}^{{d+3 \choose 3}-1}$,
the intersection $\im(\phi)\cdot H$, corresponds to the pullback
$\phi^{*}H$.

\noindent Geometrically, if $P$ is a ge\-ne\-ral point in
$\mathbb{P}^{3}$, the surfaces containing the point constitute a
hyperplane of $\mathbb{P}^{{d+3 \choose 3}-1}$. Let $T_{a_{P}}$ be
the class of the codimension 2 cycle in the Grassmannian
representing the set of lines containing the point, and $[H_{P}]$
the class of the hyperplane in $\mathbb{P}^{{d+3 \choose 3}-1}$ of
surfaces of degree $d$ containing $P$. It is clear by the way the
morphism $\phi$ has been constructed, that $\phi^{*}H_{P}=D_{P}$,
i.e. the pull back of $H_{P}$ is the divisor of rational curves
meeting the $a_{P}-$plane, or equivalently
the locus of parametrized rational ruled surfaces through $P$. 

\noindent We have seen that a morphism from $\mathbb{P}^{1}$ to
$G(2,4)$ depends on $4d+4$ parameters. If we quotient by the action of
PGL(1) on  $\mathbb{P}^{1}$,
it depends on $4d+1$ parameters. Therefore we need to
intersect $4d+1$ divisors of the first kind to fix a rational
curve on the Grassmannian, and 3 divisors of the second kind to
fix a parametrization. Now the transversality of the intersection
follows from Lemma \ref{intnum}. 
$$N_{d}=\int_{R^{0}_{d}}[D]^{4d+1}\cdot [Y]^{3}=\int_{R^{0}_{d}}
\phi^{*}[H]^{3}\cdot \phi^{*}[H_{P}]^{4d+1}=d^{3}
\int_{\mathbb{P}^{{d+3 \choose 3}-1}}[H_{P}]^{4d+1}\phi_{*}(1)$$
where 1 is the class of the total space in the Chow ring
$A(R^{0}_{d})$. The last equality is a consequence of the
projection formula. \cqd

\begin{coro}
The Severi degree of the variety of rational ruled surfaces in
$\mathbb{P}^{3}$ of degree $d$, for $d\geq 3$, in the projective
space of surfaces of degree $d$, $\mathbb{P}^{{d+3 \choose 3}-1}$,
coincides with the number of rational ruled surfaces through
$4d+1$ points.
\end{coro}

\vspace{0.5cm} In order to apply the methods of
intersection theory, a compact parameter space
is required together with some knowledge of its intersection ring.
The Quot scheme has been used many times as a smooth
compactification of the space of morphisms of a fixed degree from
$\mathbb{P}^{1}$ to a Grassmannian, \cite{Ber1}. Str\o mme in
\cite{Str} computes the Betti numbers of $R_{d}$ and gives a
description with ge\-ne\-ra\-tors and relations of its cohomology
ring. In particular, he gives a basis for its Picard group
$A^{1}(R_{d})$.

It must be checked that the solutions, i.e. the objects satisfying
all the conditions (\ref{Nd}), are in the dense open set
$R^{0}_{d}$. For this purpose the next step will be to intersect
the closures of the divisors with the boundary of $R_{d}$. The
Weil divisors defined in $R^{0}_{d}$ extend to divisors in
$R_{d}$. We shall denote by $\overline{D}_{i}$ the closures of the
divisors $D_{i}$ in $R_{d}$.

\begin{prop}
The intersection of the divisors $\overline{D}_{i}$ have a common
excess component contained in the boundary of the Quot scheme
$R_{d}$.

\end{prop}
{\it Proof.} Let us consider the universal exact sequence on
$R_{d}\times \mathbb{P}^{1}$,
$$0\rightarrow \mathcal{K}\stackrel{\varphi}{\rightarrow}
\mathcal{O}^{4}\rightarrow \mathcal{E}\rightarrow 0 \ \ \mbox{on} \ \
R_{d}\times \mathbb{P}^{1}.$$ $\mathcal{K}$ is a locally free
sheaf of rank 2, therefore $\varphi:\mathcal{K}\rightarrow
\mathcal{O}^{4}$ is a morphism of locally free sheaves. The
intersection component is given by the points $p\in R_{d}$,
$0\rightarrow N\stackrel{\varphi_{p}}{\rightarrow}
O^{4}_{\mathbb{P}^{1}}\stackrel{\bullet}{\rightarrow} E\rightarrow
0$ such that there exists $t\in \mathbb{P}^{1}$ where the morphism
$O^{4}_{\mathbb{P}^{1}}|_{t}\stackrel{\bullet}{\rightarrow}E_{t}$
is an isomorphism, i.e. the points where the rank of the map
$\varphi_{r,t}:\mathcal{K}_{r,t}\rightarrow \mathcal{O}^{4}_{r,t}$
is 0. In that case $\mathcal{E}_{\{r\}\times \mathbb{P}^{1}}$, is
not a locally free sheaf. Rather it has a non-zero torsion of
degree 2 in $t\in \mathbb{P}^{1}$, and this condition is satisfied
trivially by all the divisors $\overline{D}_{i}$. These points are
in the boundary of $R_{d}$ (see also \cite{Ber2}) and they are
parametrized by the  image of the projection $\pi:R_{d}\times
\mathbb{P}^{1}\rightarrow R_{d}$ of the determinantal variety,
$$\{(r,t)\in R_{d}\times \mathbb{P}^{1}|\ rk(\varphi_{r,t})=
0\}.$$\cqd

\begin{rem}The intersection (\ref{Nd}) does not have enumerative meaning
  on the variety $R_{d}$.  We should separate the divisors
  $\overline{D}_{i}$ outside $R^{0}_{d}$. This could be accomplished
  directly by blowing-up their intersection in $R_{d}$.  There is a
  formula for the intersection numbers in the blow ups in terms of the
  intersection rings of the base and the center, and the normal bundle
  of the center (see Co\-ro\-lla\-ry 4.2.2 and Proposition 4.1.(a),
  \cite{Ful}). The formula minimizes the amount of information needed
  to perform a single product in the intersection ring of a blow up.
  Paolo Aluffi computes in \cite{Alu} some characteristic numbers for
  smooth plane curves by using this method.
\end{rem}
The divisors $Y_{i}$ are already considered by Bertram in
\cite{Ber1} where he proves a moving lemma stating that these
varieties can be made to intersect transversely. If we denote by
$[Y]$ the associated class in the Chow ring $A(R_{d})$, then we
consider the self-intersection,
\begin{equation}
\int _{[R_{d}]}[Y]^{4d+4}\cap [R_{d}].
\end{equation}
This intersection number is computed in \cite{RRW} where it is
called the degree of the generalized Pl\"ucker embedding. It is a
Gromov-Witten invariant and can also be obtained  by computing in the
Quot scheme by means of the formulas of Vafa and Intriligator,
\cite{Ber1}.

\section{Interpretation of  $ N_{d}$ as a Gromov-Witten in\-va\-riant}

\setcounter{equation}{0} The Kontsevich moduli space
$\overline{M}_{0,n}(G(2,4),d)$ parametrizes isomorphism classes of
stable maps $f$ from nodal, $n-$pointed trees of
$\mathbb{P}^{1}$'s $(C,p_{1},\ldots p_{n})$ into $G(2,4)$ of
Pl\"{u}cker degree $d$, where
$f_{*}([C])=dT_{3}$. The stability condition
implies that any $\mathbb{P}^{1}$ contracted by $f$ has at least
three nodes or marked points. This space has dimension $4d+n+1$
and provides a compactification of the space of morphisms from
$\mathbb{P}^{1}$ to $G(2,4)$.

\noindent The numbers $N_{d}$ occur as intersection numbers on the
space \footnotesize{$\overline{M}_{0,4d+4}(G(2,4),d)$}.
\normalsize For each marked point $1\leq i\leq 4d+4$, there is a
canonical evaluation map:
\begin{equation}
\pi_{i}:\overline{M}_{0,4d+4}(G(2,4),d)\rightarrow G(2,4) \ \ \
\pi_{i}([C,p_{1}\ldots ,p_{n},f])=f(p_{i}).
\end{equation}
A product in the ring $A^{*}(\overline{M}_{0,n}(G(2,4),d))$ is
given by
\begin{equation}
\pi_{1}^{*}(\gamma_{1})\cup \ldots \cup \pi_{n}^{*}(\gamma_{n})\in
A^{*}(\overline{M}_{0,n}(G(2,4),d)),
\end{equation}
where $\gamma_{i}$ are cycles on $G(2,4)$.

\noindent If $\sum
\codim(\gamma_{i})=\dim(\overline{M}_{0,n}(G(2,4),d))$, {\it
the Gromov-Witten invariant}

\noindent $I_{0,n,d}(\gamma_{1},\ldots,\gamma_{n})$ is defined as the top
degree class:
\begin{equation}
I_{0,n,d}(\gamma_{1},\ldots,\gamma_{n})=\int_{[\overline{M}_{0,n}(G(2,4),d)]^{vir}}
\pi_{1}^{*}(\gamma_{1})\cup \ldots \cup \pi_{n}^{*}(\gamma_{n}).
\end{equation}

\noindent If $\gamma_{i},\beta_{j}\in H^{*}(G(2,4),\mathbb{Z})$
are cohomology classes satisfying:
\begin{equation}\label{dim}\sum_{i=0}^{k} \deg\,\gamma_{i}+
\sum_{j=0}^{l} \deg\,\beta_{j}=4d+4+l \end{equation} then {\it the
mixed invariant}
$\varphi_{0,4d+4,d}\,(\gamma_{1},\ldots,\gamma_{k}|\beta_{1},\ldots
\beta_{l})$  is defined as the number of maps
$\mu:\mathbb{P}^{1}\rightarrow G(2,4)$ in
$\overline{M}_{0,4d+4}(G(2,4),d)$ such that
$\mu_{*}[\mathbb{P}^{1}]=d T_3$, $\mu(p_{i})$ lies in a cycle dual to
$\gamma_{i}$ and $\mu(\mathbb{P}^{1})$ meets a cycle dual to
$\beta_{j}$.

\vspace{0.5cm}  As in the previous section we fix the action of
the group PGL(3) on  $G(2,4)$ and we choose general $a_{i}-$planes
for $i=1,\ldots 4d+1$ and 3 general hyperplanes $j=1,2,3$ in
$G(2,4)$ with $[T_{a_{i}}]$ and $[T_{1_{j}}]$ their cohomology
classes respectively in the Chow ring $A(G(2,4))$.

\begin{theorem}
The Gromov-Witten invariant,
$$I_{0,4d+4,d}(T_{1},T_{1},T_{1},T_{a_{1}},\ldots,T_{a_{4d+1}}),$$
coincides with the intersection number $N_{d}$ ($d \geq 3$), that is, with the
number of rational curves in $G(2,4)$ of degree $d$ which have
non-empty intersection with all $a_{i}-$planes $(i=1,\ldots
,4d+1)$, $j-$hyperplanes $(j=1,2,3)$.
\end{theorem}
{\it Proof.} Let $Y$ be the subscheme in
$\overline{M}_{0,4d+4}(G(2,4),d) $ defined as the intersection
product,
\begin{equation}\label{product}
Y=\pi _{4d+2}^{-1}(T_{3})\cap \pi_{4d+3}^{-1}(T_{1})\cap\pi
_{4d+4}^{-1}(T_{1})\cap \pi_{1}^{-1}(T_{a_{1}})\cap \ldots \cap
\pi_{4d+1}^{-1}(T_{a_{4d+1}}).
\end{equation}

\noindent The evaluation maps are flat by generic flatness
\cite{Alt}, therefore the inverse image of a hyperplane is of
codimension 1 and the inverse image of an $a$-plane is of
codimension 2. We have seen that
$\dim(\overline{M}_{0,4d+4}(G(2,4),d))=8d+5$ and that the
codimension of the total intersection by the previous observation
is $2\cdot (4d+1)+3=8d+5$. \noindent $Y$ is the locus of maps
$(C;P_{1},\ldots,P_{4d+4};\mu)$ such that $\mu(P_{i})\in
T_{a_{i}}$ for $i=1\ldots 4d+1$ and $\mu(P_{4d+1+j})\in T_{1_{j}}$
for $j=1,2,3$. Since $G(2,4)$ is a homogeneous variety, the number
of these maps equals
$I_{0,4d+4,d}(T_{1},T_{1},T_{1},T_{a_{1}},\ldots,T_{a_{4d+1}})$,
(see Lemma 7.14 of [FP]). Since the number of marked points $n$ is
$\geq 3$, we can work in the variety of morphisms $R^{0}_{d}$ (see
\ref{finemod}). If we fix three of the marked points to be
$0,1,\infty \in \mathbb{P}^{1}$ and let $V=\mathbb{P}^{1}\times
\stackrel{4d+1}{\ldots}\times \mathbb{P}^{1}\backslash
\{\triangle_{i,j},L_{0,i},L_{1,i},L_{\infty,i}\}$ where
$\triangle_{i,j}$ is the large diagonal determined by factors $i$
and $j$, $L_{0,i}$ is the locus where the $i^{th}$ factor is $0\in
\mathbb{P}^{1}$ and $L_{1,i},L_{\infty,i}$ are defined similarly,
there is a universal family of Kontsevich stable degree $d$ maps
of $4d+4-$pointed automorphism free curves with smooth and
irreducible source,
\begin{equation}
\mathbb{P}^{1}\times V\times R^{00}_{d}\rightarrow G(2,4).
\end{equation}
By the universal property there is an injection $$V\times
R^{00}_{d}\hookrightarrow \overline{M}_{0,4d+4}(G(2,4),d).$$ Now
arguing as in the proof of Lemma \ref{intnum}, it follows by
Kleiman's Theorem that the intersection (\ref{product}) is a
finite set of reduced points and is supported in the open set
$V\times R^{00}_{d}$.

\noindent Next we prove that counting stable maps is the same as
counting rational curves. First we check there is no repetition,
that is, the solution curves intersect any cycle in just one
point. In the case of $T_{1_{j}}$, this type of repetition is
unavoidable. By Bezout's theorem, a curve of degree $d$ will
always meet a codimension one space. This is just the
divisor-property of the Gromov-Witten invariants. This gives us a
contribution factor $d^{3}$,
\begin{equation}\label{divisor}
I_{0,4d+4,d}(T_{1},T_{1},T_{1},T_{a_{1}},\ldots,T_{a_{4d+1}})=d^{3}I_{0,4d+1,d}
(T_{a_{1}},\ldots,T_{a_{4d+1}}).
\end{equation}
\noindent For the cycles of codimension 2, another transversality
argument (with res\-pect to the PGL(3)-action) implies each
solution map intersects the cycle in only one point. This part is
similar to the analogous result for plane curves (see Lemma 3.5.5.
of [KV]).

\noindent The other behavior we want to exclude is when the same
curve passes twice through the same point. In order to prove this,
we have to see that the locus of maps $\mu$ for which
$\mu^{-1}(\mu(P_{i}))$ contains at least one point distinct from
$P_{i}$, or equivalently, the locus of maps with nodes, is of
positive codimension. Let us consider the diagram,
\begin{equation}\label{diag}
\begin{array}{ccc}
T & \rightarrow      & \triangle \\
\downarrow      &  & \downarrow \\
X:=R^{0}_{d}\times\mathbb{P}^{1}\times\mathbb{P}^{1} &
\stackrel{e}{ \rightarrow} & G(2,4)\times G(2,4)
 \end{array} \end{equation}
where $\triangle$ is the diagonal in $G(2,4)\times G(2,4)$, and
$T$ is the fiber product $X \times_{G\times G} \triangle$ which is
the locus of morphisms with at least one node. We observe that the
action of the group PGL(3)$\times$PGL(3) on $G(2,4)\times G(2,4)$
is not transitive, in fact there are two orbits, a dense orbit and
the diagonal orbit. A slightly more general result on the
transversality concerning the orbits of an action \cite{Spe} apply
to see that $T$ has positive codimension.
\cqd

\vspace{0.5cm} We define $Q_{d}$ to be the Gromov-Witten invariant
$I_{0,4d+1,d}(T_{a_{1}},\ldots,T_{a_{4d+1}})$.
\begin{coro} $Q_{d}$ counts the number of rational ruled surfaces
in $\mathbb{P}^{3}$ through $4d+1$ points for each degree $d\geq
3$.
\end{coro}

\begin{rem} We are counting maps of degree
$d$, $f:\mathbb{P}^{1}\rightarrow G(2,4)$ verifying
$f(\mathbb{P}^{1})\cap T_{a_{i}}\neq \emptyset \ \
i=1,\ldots,4d+1$ and $f(P_{j})\in T_{1_{j}}\ \ j=1,2,3$ and the
cycles verify the assumption (\ref{dim}) on the dimension, so
$N_{d}$ is also the mixed invariant
$$\varphi_{0,4d+4,d}(T_{3},T_{3},T_{3}|T_{a_{1}},\ldots,T_{a_{4d+1}}).$$
\end{rem}

\section{Effective computation of degree and quantum-cohomology}

We now describe the quantum cohomology ring of $G(2,4)$,
$QH^{*}(G(2,4))$. The basic reference for this section is
\cite{FI}. We introduce the new variable $T_{0}$ for the class of
the total G(2,4).

The classical ring with unit $T_{0}$ is given by the relations:
$$T_{1}^{2}=T_{a}+T_{b}$$
$$T_{1}T_{a}=T_{1}T_{b}=T_{3}$$
$$T_{1}T_{3}=T_{a}^{2}=T_{b}^{2}=T_{4}.$$

The Gromov-Witten invariants count the numbers:
\begin{equation}\label{identification}
N(\alpha,\beta,\gamma,\delta;d)=\,
I_{0,4d+1,d}(T_{a}^{\alpha},T_{b}^{\beta},T_{3}^{\gamma},T_{4}^{\delta}),
\end{equation}
which are nonzero when
$\alpha+\beta+2\,\gamma+3\,\delta=dim\,\overline{M}_{0,0}(G(2,4),d)=4d+1,\
\ \rm with\ \alpha, \beta, \gamma, \delta \geq 0$.

We define the numbers $g_{ij}=\int_{G}T_{i}\cup T_{j},\ \ \ i,j\in
\{1,a,b,3,4\}$, and $(g^{ij})$ for the inverse matrix. The only
non-vanishing elements of the intersection form are
$(g_{04}=g_{aa}=g_{bb}=g_{13}=1)$.
\begin{equation}
T_{i}\cup T_{j}=\sum_{e,f}\int_{G}(T_{i}\cup T_{j}\cup T_{e})\,
g^{ef}\, T_{f}=\sum_{e,f}\,I_{0}(T_{i}T_{j}T_{e})\,g^{ef}T_{f}
\end{equation}
A  \lq\lq quantum deformation" of the cup multiplication is
defined by allowing nonzero classes.

 Now we introduce the {\bf deformation
parameters} $y_{0},y_{1},y_{a},y_{b},y_{3}$, $y_{4}$ and the genus
0 free energy $F$ splits into $F=f_{cl}+f$, where the classical
potential function is given by
\begin{equation} \label{potencial}
f_{cl}(y_{0},y_{1},y_{a},y_{b},y_{3},y_{4})=\frac{1}{2}\,y_{0}\,(y_{4}y_{0}+
y_{a}^{2}+y^{2}_{b})+\frac{1}{2}y_{1}^{2}(y_{a}+y_{b})+y_{0}y_{1}y_{3},
\end{equation} and the quantum potential by
\begin{equation}\label{q-potencial}
f=\sum_{\alpha, \beta, \gamma, \delta >0 \
\alpha+\beta+2\gamma+3\delta=4d+1}N(\alpha,\beta,\gamma,\delta;d)
\frac{y_{a}^{\alpha}}{\alpha!}\frac{y_{b}^{\beta}}{\beta!}\frac{y_{3}^{\gamma}}{\gamma!}
\frac{y_{4}^{\delta}}{\delta!}\,e^{dy_{1}}.
\end{equation}

Quantum corrections are due to non-trivial maps
$\mathbb{P}^{1}\rightarrow G(2,4)$ with the necessary markings to
rigidify them. The splitting of the free energy into $F_{cl}$ and
$F$ is according to maps for which the image of $\mathbb{P}^{1}$
is a point or an irreducible curve.
\begin{equation}
F_{ijk}:=\frac{\partial^{3}\,F}{\partial\,y_{i}\partial\,y_{j}\partial\,y_{k}},\
\ \ 0\leq i,j,k\leq m
\end{equation}
The {\bf quantum product} is defined as:

\begin{equation}
T_{i}*\,T_{j}=\sum_{e,f}F_{ijk}g\,^{ef}T_{f}\ \mbox{where} \  T_{f} \mbox{ is
Poincar\'{e} dual to} \  T_{e}
\end{equation}
This product is extended $Q[[y]]$-linearly to the $Q[[y]]$-module
$A^{*}G\otimes_{\mathbb{Z}}Q[[y]]$ thus making it a
$Q[[y]]$-algebra. The product is commutative, since the partial
derivatives are symmetric in the subscripts. It is easy to see
that $t_{0}=1$ is a unit.The essential point is the
associativity, \cite{FP}, which allows us to obtain interesting
relations among the quantum-numbers
$N(\alpha,\beta,\gamma,\delta;d)$.

There are 55 relations expressing the associativity of the
deformed ring of the Grassmannian and the obvious
symmetry relations
$N(\alpha,\beta,\gamma,\delta|d)=N(\beta,\alpha,\gamma,\delta|d)$.
Via the identification (\ref{identification}), the numbers $Q_{d}$
are of the form $N(4d+1,0,0,0|d)$. To compute these numbers we
have applied Farsta, a program due to Andrew Kresch
(http://www.maths.warwick.ac.uk/~kresch/).

The following table lists the number of rational ruled surfaces in
$\mathbb{P}^{3}$ through $4d+1$ points for each degree $3\leq d\leq
9$. In the case $d=2$, $Q_{2}$ is twice the number of quadrics through
9 points, that is 1, since as we saw in Section \ref{themorph} a
quadric has 2 different rulings.

\vspace{0.5cm}
\begin{center}
\begin{tabular}{|c|r|r|}
\hline $d$ & $Q_{d}$  & $4d+1$   \\
\hline
 1 & 0& 5 \\
 2 & 2& 9 \\
 3 & 504& 13 \\
 4 & 1044120 & 17 \\
 5 & 5335687360 & 21 \\
 6 & 67992124121040 & 25 \\

 7 & 1743784747544391896 & 29 \\
 8 & 82475300124495938244352 & 33 \\
 9 & 6608238869716397977928547520 & 37 \\
 \hline
 \end{tabular}
 \end{center}
 \vspace{0.5cm}
The case $d=3$ has been studied by Vainsencher and Coray in
\cite{CV} where they showed there are 504 rational ruled cubic
surfaces passing through 13 general points or, equivalently, 504
rational curves in the Grassmannian $G(2,4)$ meeting 13 general
$a-$planes.

\end{document}